\documentclass[a4paper,12pt]{article}

\title{\textbf{CLASS FIELD THEORY AND ARITHMETIC OF ABELIAN VARIETIES OVER LOCAL FIELDS}}
\author{Christopher Stephen Hall}
\date{}
\usepackage{amsmath,amssymb,amsfonts} % Typical maths resource packages
\usepackage{graphics}                 % Packages to allow inclusion of graphics
\usepackage{color}                    % For creating coloured text and background
\usepackage{hyperref} 
\usepackage[all, cmtip]{xy}
\usepackage{tipa}
\usepackage{geometry}
\begin{document}

\maketitle

{\small \noindent ABSTRACT. We use knowledge of local fields to adapt Jonathan Lubin and Michael Rosen's proof of Mazur's Proposition $4.39$. This changes the result about abelian varieties from only working over local fields with a finite residue field to working with local fields with an arbitrary perfect residue field of positive characteristic. We then briefly discuss the Local Class Field Theory implications of such information} 
\\
\\
\\

\section{Introduction}

This paper is about class field theory for abelian varieties with ordinary good reduction over local fields. It is related to and further extends Jonathan Lubin and Michael I. Rosen's 1977 paper ``The Norm Map for Ordinary Abelian Varieties'' \cite{Rubin1}. It will endeavour to generalise the main results of that paper, one of which is a reproof of Mazur's Proposition $4.39$ \cite{Mazur1}, to a wider collection of complete discrete valuation fields.

In the original paper Lubin and Michael did not reference Local Class Field Theory despite the fact that the topic involves Galois extensions of Local Fields, and one of the results they get can be seen as an analogue to a major result of Local Class Field Theory. I, however, shall be burrowing ideas from the topic, in particular from Professors Fesenko and Vostokov's book ``Local Fields and their Extensions''; in which chapter $5$ has a brief introduction to Local Class Field Theory when the Local Fields has an arbitrary perfect residue field of positive characteristic \cite{fesenko1}. I will also be utilising Fesenko's paper ``Local Class Field Theory: Perfect Residue Field Case'' as this goes into that topic in a lot more detail than his book does \cite{Ivan1}, in fact the section on the topic in the book explicitly states how that paper should be used to explore the subject matter further \cite{fesenko1}.

Before I get to the main result I should explain what we are going to work with.

We are letting $K$ be a complete discrete valuation field of characteristic $0$ and that $\mathtt{char}(k) = p > 0$, where $k = \overline{K}$is a perfect field.

We want to prove the following result:

\bigskip

\textbf{Theorem 0:} \textit{Let $A$ be a d-dimensional abelian variety with good reduction over $K$ with family of twist matrices $u_{J}$, the family shall be explained later. If $L/K$ is a totally ramified $\mathbb{Z}_{p}$-extension of $K$ with $\mathcal{N}_{L/K}(A(L)) = \bigcap N_{K_{n}/K}(A(K_{n}))$, where $L/K_{n}/K$ and $[K_{n} : K] = p^{n}$, then the following exact sequence can be constructed:}

\[ \oplus_{j \in J} \mathbb{Z}_{p}^{d}/(I - u_{J})(\oplus_{j \in J} \mathbb{Z}_{p}^{d}) \to A(K)/\mathcal{N}_{L/K}(A(L)) \to A(k)_{p} \to (0) \]

\noindent\textit{Here $A(k)_{p}$ is the group of $p$-torsion of $A(k)$, so the elements of $A(k)$ whose order is a finite power of $p$. $A(k)$, meanwhile, is the group of the $k$-value points of the reduced abelian variety $\overline{A}$, which is defined over $k$.}

\bigskip

Other than the first group this sequence is a copy of the exact sequence in Lubin and Rosen's paper \cite{Rubin1}; the only difference being the first term. I shall explain what $J$ is in the next section but it is related to the fact that I am exploring this topic from the perspective of Local Class Field Theory. After proving the above sequence is exact I will showcase that the similarity is not a coincidence and that it is a generalisation of the exact sequence that appears in Lubin and Rosen's paper, by that I mean that if you assume that $k$ is finite the above sequence simplifies to the one Lubin and Rosen described.

The above exact sequence is just the beginning of what one can explore in this topic, in particular when one looks at the subject from the point of view of Local Class Field Theory. For that reason the last section of this paper shall be given over to briefly considering the further ways one could take the topic of a Local Class Field Theory approach to abelian varieties over Local Fields with a perfect residue field.

Finally, despite the fact that I only cite $4$ other mathematical works in this paper it may be clear that I have only talked about $3$ of them in the above introduction. This is because the third one, which is in fact Mazur's paper ``Rational Points of Abelian Varieties with Values in Towers of Number Fields'' in which Lubin and Rosen got the original exact sequence from \cite{Mazur1}, is only utilised for a single result regarding the nature of abelian varieties with good reduction over Local Fields and other than that is irrelevant to the mathematics I am utilising here.

\bigskip

\section{Notation}

Before we can start with the proof of the above theorem we first must define some of the notation that we will be using. As this is an expansion of ``The Norm Map for Ordinary Abelian Varieties'' I will be using the same notation that Lubin and Rosen uses \cite{Rubin1}, for the sake of continuity and for the ease of the reader.

Let $E/K$ be either an arbitrary algebraic extension of $K$ or the completion of one. Either way $E$ is a discrete valuation ring and we will denote its ring of integers by $\mathcal{O}(E)$. We will also let $\mathcal{U}(E)$ and $\mathcal{U}^{1}(E)$ be the units and principle units of $E$ respectively.

Now, we will let $T/K$  be the completion of the maximal unramified $p$-extension of $K$, with $\overline{T}$ being the residue field of $T$. As $k$ has characteristic $p$ we get that $\mathtt{Gal}(T/K)$ is a free $\mathbb{Z}_{p}$-group. Let $\phi_{j}, j \in J$ be a collection of topological generators of $\mathtt{Gal}(T/K)$.

Next, we have the totally ramified $p$-extension $L/K$, but rather than working directly with $\mathtt{Gal}(L/K)$ we instead work with the group of continuous homomorphisms $\mathtt{Gal}(L/K)\,\widehat{•}$. If $G$ is a $p$-group then $G\,\widehat{•} = \mathtt{Hom}_{\mathbb{Z}_{p}}(\mathtt{Gal}(T/K), G)$. This definition is taken from chapter $5$ of ``Local Fields and their Extensions'' \cite{fesenko1}. We also have that if $G$ is abelian then $G\,\widehat{•} \cong \oplus_{j \in J} G$.

Let $F$ be a $d$-dimensional toroidal formal group over $\mathcal{O}(K)$ and let $f: F \to \hat{\mathbb{G}}_{m}^{d}$ be an isomorphism over $\mathcal{O}(K^{\mathtt{ur}})$, where $K^{\mathtt{ur}}$ is the completion of the maximal unramified extension of $K$. Such an isomorphism exists from the definition of a toroidal formal group.

$f$ can be described by a power series and by applying $\phi_{j}$, for some $j \in J$, to the coefficients of the power series we end up with the isomorphism of formal groups $f^{\phi_{j}}: F \to \hat{\mathbb{G}}_{m}^{d}$. Let $u_{j}$ be the invertible $d \times d$ matrix over $\mathbb{Z}_{p}$ that corresponds to the automorphism of $\hat{\mathbb{G}}_{m}^{d}$, $f^{\phi_{j}} \circ f^{-1}$. The $u_{j}$ is called a twist matrix of $F$ for $j \in J$. We denote such a family of twist matrices by $u_{J}$. As an aside, this explains the notation that ``$u_{J}$'' that I wrote in the exact sequence in the introduction.

Finally, if $L/K$ is a totally ramified Galois extension then we may extend any of the $\phi_{j}$'s, for $j \in J$, to $LT$ by having $\phi_{j}$ act trivially on elements of $L$.

\bigskip

\subsection{The Module $V(L)$}

Let $L/T$ be a finite totally ramified Galois extension. Extend $\phi_{j}, j \in J$ to $LT$; then we may define the group:

\[ V(L) = V_{u_{J}}(L) = \{\alpha \in \mathcal{U}^{1}(LT)^{d} : \alpha^{\phi_{j}} = \alpha^{u_{j}}, \forall j \in J \} \]

For the sake of notation, we shall denote $V_{u_{J}}(L)$ as $V(L)$ when $u_{J}$ is obvious.

We have that the $\phi_{j}$'s act diagonally on $\alpha$ and the $u_{j}$'s act in the obvious way. The actions of the $\phi_{j}$'s commute with the elements of $\mathtt{Gal}(L/K)$ and thus $V(L)$ is a $\mathtt{Gal}(L/K)$-module. What is more $K(\mathcal{O}(L)) \cong V(L)$ as $\mathtt{Gal}(L/K)$-modules.

To show this let $f: F \to \hat{\mathbb{G}}_{m}^{d}$ be an isomorphism over $\mathcal{O}(T)$ such that $u_{j} = f^{\phi_{j}} \circ f^{-1}$ for all $j \in J$, we know such an $f$ exists from the definition of $u_{j}$. Now, let $\alpha$ be in $F(\mathcal{O}(LT))$ and suppose $f(\alpha) \in V(L)$; $f$ defines an isomorphism of groups $F(\mathcal{O}(LT))$ and $\mathcal{U}^{1}(LT)^{d}$. We have that for all $j \in J$:

\[f(\alpha)^{\phi_{j}} = f(\alpha)^{u_{j}} = f^{\phi_{j}} \circ f^{-1} \circ f(\alpha) = f^{\phi_{j}}(\alpha)\] 

\noindent This means that $f(\alpha)^{\phi_{j}} = f^{\phi_{j}}(\alpha^{\phi_{j}})$ and since $f^{\phi_{j}}$ is an isomorphism we have that $f^{\phi_{j}}(\alpha^{\phi_{j}}) = f^{\phi_{j}}(\alpha)$ implies that $\alpha^{\phi_{j}} = \alpha$. So $f(\alpha) \in V(L)$ tells us that $\alpha^{\phi_{j}} = \alpha$ for all $j \in J$. To finish this part of the proof we note that $\alpha \in F(\mathcal{O}(L))$ if and only if $\alpha^{\phi_{j}} = \alpha$ for all $j \in J$. We, therefore, have that if $f(\alpha) \in V(L)$ then $\alpha \in F(\mathcal{O}(L))$. The result that $f(\alpha) \in V(L)$ follows from $\alpha \in F(\mathcal{O}(L))$ is proved in a similar fashion.

I will note that Mazur's ``Rational Points of Abelian Varieties with Values in Towers of Number Fields'' states that since $k$ is perfect then if $A$ is a $d$-dimensional abelian variety over $K$ with good reduction then the formal group, $\hat{A}$, corresponding to $A$ is toroidal \cite{Mazur1}.

\bigskip

\section{Theorem 1}

\textit{There is an isomorphism:}

\[ V_{u_{J}}(K)/\mathcal{N}_{L/K}(V_{u_{J}}(L)) \cong ((G^\mathtt{ab})^{d})\,\widehat{•}\,/(I - u_{J})(((G^\mathtt{ab}\,)^{d})\,\widehat{•}\,) \]

\noindent \textit{where $L/K$ is a finite totally ramified $p$-extension, $G = \mathtt{Gal}(L/K)$ and $u_{J}$ is the family of twist matrices relating to $\hat{A}$.}

\bigskip

\textbf{Proof:} Proving the above isomorphism will be a major stepping stone towards showing my main result. Also, readers with a copy of Lubin and Rosen's paper will notice that this is very similar to Theorem $1$ of that paper.

If we let $E = LT$ then we get that $\mathtt{Gal}(E/T) \cong G$ and $\mathtt{Gal}(E/L) \cong \mathtt{Gal}(T/K)$, with $\mathtt{Gal}(E/L)$ being topologically generated by the modified versions of the $\phi_{j}$'s that we saw when discussing $V(L)$. Finally, it should be noted that, $E/L$ is the completion of the maximal unramified $p$-extension of $L$.

To continue the proof of Theorem $1$ we first need to deal with a few lemmas.

\bigskip

\section{Lemma 1}

\textit{The $G$-module $E^{\times}$ is cohomologically trivial.}

\bigskip

\textbf{Proof:} Firstly, as $E$ and $T$ have algebraically $p$-closed residue fields then the norm map $N_{E/T}$ is surjective \cite{fesenko1}. This means that for any subgroup $H \subseteq G$, we have that $\hat{H}^{0}(H, E^{\times}) = (0)$. Likewise, by Hilbert's Theorem $90$, $H^{1}(H, E^{\times}) = (0)$. So, the cohomology groups $H^{n}(H, E^{\times})$, for any subgroup $H$ of $G$, vanish in the successive dimensions $n= 0$ and $n = 1$; thus $E^{\times}$ is a chomologically trivial $G$-module.

\bigskip

\section{Lemma 2}

\textit{There is an exact sequence:}

\[ (1) \to (G^{\mathtt{ab}})\,\widehat{•} \to (\mathcal{U}^{1}(E)/(I_{G}\mathcal{U}^{1}(E)))\,\widehat{•} \to (\mathcal{U}^{1}(T))\,\widehat{•} \to (1) \]

\bigskip

\textbf{Proof:} Since $G$ is a $p$-group, and therefore so is $G^{ab}$, we have that $(G^{\mathtt{ab}})\,\widehat{•} \cong \oplus_{j \in J} G^{\mathtt{ab}}$. The same holds for the $p$-groups $\mathcal{U}^{1}(T)$ and $\mathcal{U}^{1}(E)/(I_{G}\mathcal{U}^{1}(E))$. Here $I_{G}$ is the kernel of the augmentation map from $\mathbb{Z}[G]$ to $\mathbb{Z}$ that sends $\sum_{g \in G} n_{g}g \in \mathbb{Z}[G]$, with $n_{g} \in \mathbb{Z}$, to $\sum_{g \in G} n_{g} \in \mathbb{Z}$.

The next thing to do is to show that for any $\phi_{j}$, $j \in J$, that we have the exact sequence:

\[ (1) \to G^{\mathtt{ab}} \to \mathcal{U}^{1}(E)/(I_{G}\mathcal{U}^{1}(E)) \to \mathcal{U}^{1}(T) \to (1) \]

To begin with, fix an $j \in J$ to work with. From the previous lemma and the exact sequence:

\[ (1) \to \mathcal{U}(E) \to E^{\times} \to \mathbb{Z} \to (1) \]

\noindent we have the isomorphisms:

\[ \mathbb{G}^{\mathtt{ab}} \cong H^{-2}(G, \mathbb{Z}) \cong H^{-1}(G, \mathbb{Z}) \cong _{\mathcal{N}}\mathcal{U}(E)/(I_{G}\mathcal{U}(E)) \]

From the surjectivity of the norm maps, $N_{E/T}(\mathcal{U}(E)) = \mathcal{U}(T)$. Thus we get the exact sequence:

\[ (1) \to G^{\mathtt{ab}} \to \mathcal{U}(E)/(I_{G}\mathcal{U}(E)) \to \mathcal{U}(T) \to (1) \]

\noindent This is nearly the exact sequence that we are after.

$\overline{T}$ is the separable $p$-closure of $k = \overline{K}$. As $E = LT$ and $L/K$ is totally ramified, $\overline{E} = \overline{T}$. $\overline{E}$ has characteristic $p$ and thus is uniquely $p$-divisible, so for all $m \in \mathbb{Z}$, $H^{m}(G, \overline{E}^{\times}) = 0$. This, when combined with the exact sequence:

\[(1) \to \mathcal{U}^{1}(E) \to \mathcal{U}(E) \to \overline{E}^{\times} \to (1)\] 

\noindent means that for all $m \in \mathbb{Z}$, we have $\hat{H}^{m}(G, \mathcal{U}^{1}(E)) \cong \hat{H}^{m}(G, \mathcal{U}(E))$.

Set $m = -1$ to get $\mathcal{U}(E)/(I_{G}\mathcal{U}(E)) \cong  \mathcal{U}^{1}(E)/(I_{G}\mathcal{U}^{1}(E))$. Then, if we let $m = 0$, we observe that $\hat{H}^{0}(G, \mathcal{U}^{1}(E)) \cong \hat{H}^{0}(G, \mathcal{U}(E))$. The latter group is isomorphic to $\mathcal{U}(E)^{G}/N_{E/T}(\mathcal{U}(E)) = (0)$, from the surjectivity of the norm map and the fact that the elements of $\mathcal{U}(E)$ that $G$ acts trivially on are precisely those in $\mathcal{U}(T)$. This, alongside with $\mathcal{U}^{1}(E)^{G} = \mathcal{U}^{1}(T)$ and $N_{E/T}(\mathcal{U}^{1}(E)) \subseteq \mathcal{U}^{1}(T)$, shows that the norm map $N_{E/T}$ maps $\mathcal{U}^{1}(E)$ surjectively onto $\mathcal{U}^{1}(T)$.

Combining the previous results with the intermediate exact sequence we got before gives us the exact sequence:

\[ (1) \to G^{\mathtt{ab}} \to \mathcal{U}^{1}(E)/(I_{G}\mathcal{U}^{1}(E)) \to \mathcal{U}^{1}(T) \to (1) \]

\noindent which is what we were trying to prove.

However, we have only worked with $\phi_{j}$ for an arbitrary $j \in J$. We can change that by acting on a term by term basis for every $j \in J$ and thus end up with the composite exact sequence:

\[ (1) \to (G^{\mathtt{ab}})\,\widehat{•} \to (\mathcal{U}^{1}(E)/(I_{G}\mathcal{U}^{1}(E)))\,\widehat{•} \to (\mathcal{U}^{1}(T))\,\widehat{•} \to (1) \]

\noindent that we were after.

\bigskip

\section{Lemma 3}

\textit{The following sequence:}

\[ \xymatrix{ (1)  \ar[r] & V(K) \ar[r] & (\mathcal{U}^{1}(T)^{d})\,\widehat{•} \ar[r]^{\phi_{J} - u_{J}} & (\mathcal{U}^{1}(T)^{d})\,\widehat{•} \ar[r] & (1) } \]

\noindent \textit{is exact.}

\textit{The map $\phi_{J} - u_{J}: (\mathcal{U}^{1}(T)^{d})\,\widehat{•} \to (\mathcal{U}^{1}(T)^{d})\,\widehat{•}$ is defined to be the homomorphism that applies $\phi_{j} - u_{j}$ to the $j$'th summand of $\oplus_{j \in J} \mathcal{U}^{1}(T) = \mathcal{U}^{1}(T)^{d}$.}

\bigskip

\textbf{Proof:} The only part of the proof that is not obvious is showing that:

\[ \xymatrix{ (\mathcal{U}^{1}(T)^{d})\,\widehat{•} \ar[r]^{\phi_{J} - u_{J}} & (\mathcal{U}^{1}(T)^{d})\,\widehat{•} } \]

\noindent is a surjective map.

Since $\mathcal{U}^{1}(T)^{d}$ is abelian, $(\mathcal{U}^{1}(T)^{d})\,\widehat{•} \cong \oplus_{j \in J} \mathcal{U}^{1}(T)^{d}$. And as $\phi_{J} - u_{J}$ acts on each summand of $\oplus_{j \in J} \mathcal{U}^{1}(T)^{d}$ independently we may fix an arbitrary $j \in J$ and instead prove the surjectivity of $\phi_{j} - u_{j}: \mathcal{U}^{1}(T)^{d} \to \mathcal{U}^{1}(T)^{d}$.

For ease of notation we shall from now on write $\phi_{j}$ as $\phi$ and $u_{j}$ as $u$. To make matters even easier the  filtration of $\mathcal{U}^{1}(T)^{d}$, which is $(\mathcal{U}^{n}(T))^{d}$ for $n \ge 1$, is preserved by $\phi - u$ and thus we will only need to check is that the induced map $\overline{\phi} - \overline{u}: \overline{T}^{d} \to \overline{T}^{d}$ is surjective.

\bigskip

Note that as the proof is the same, doing so will also prove that the sequence:

\[ \xymatrix{ (1)  \ar[r] & V(L) \ar[r] & (\mathcal{U}^{1}(E)^{d})\,\widehat{•} \ar[r]^{\phi_{J} - u_{J}} & (\mathcal{U}^{1}(E)^{d})\,\widehat{•} \ar[r] & (1) } \]

\noindent is exact as well.

\bigskip

For the next part of the proof of Lemma $3$ we shall be dealing with two separate cases, based on the value of $d$.

\bigskip

\subsection{When d equals 1}

To begin with, we shall deal with the case when $d = 1$; thus $F$ is a one-dimensional formal group over $\mathcal{O}(K)$. This means that we will be proving the surjectivity of $\overline{\phi} - \overline{u}: \overline{T} \to \overline{T}$.

As $d = 1$, $u$ is an invertible element of $\mathbb{Z}_{p}$ and thus $\overline{u} \in \mathbb{F}_{p}^{\times}$. The case where $\overline{u} = 1$, and thus $u = 1$, is covered in chapter $5$ of ``Local Fields and their Extension'' \cite{fesenko1}, so instead we will assume that $\overline{u} \not= 1$.

Let $\alpha \in \overline{T}$. Now $T/K$ is a Galois extension and therefore so is $\overline{T}/k$; thus there is a finite subfield $M/k$ of $\overline{T}/k$ such that $\alpha \in M$. $\overline{T}/k$ is the maximal $p$-extension of $k$ and thus $M/k$ has degree $p^{n}$ for some non-negative integer $n$.

The first case to consider is if $\overline{\phi}(\alpha) = \alpha$. Then $(\overline{\phi} - \overline{u})(\alpha) = (1 - \overline{u})(\alpha)$. Since $\overline{u} \not= 1$ we may set $\beta = (1 - \overline{u})^{-1}(\alpha) \in M \subseteq \overline{T}$. Finally, we observe that $(\overline{\phi} - \overline{u})(\beta) = \alpha$.

If $\overline{\phi}$ does not act trivially on $\alpha$, then for ease of later mathematics we shall multiply the map $( \overline{\phi} - \overline{u})$ by the constant $-\overline{u}^{-1}$ to get the new map $(1 - \overline{u}^{-1}\overline{\phi})$.

In this case what we should consider is that as $\alpha \in M$ then $\overline{\phi}^{p^{n}}(\alpha) = \alpha$. Now, for all non-negative integers $n$, we have $\overline{u}^{p^{-n}} = \overline{u} \not= 1$.

Set $\gamma \in M$ to equal the sum $ (1 + \sum_{m = 1}^{m = p^{n} - 1} \overline{u}^{-m} \overline{\phi}^{m})(\alpha)$. By telescoping, $(1 - \overline{u}^{-1}\overline{\phi})(\gamma) = (1 - \overline{u}^{p^{-n}})(\alpha)$. We have $(1 - \overline{u}^{p^{-n}}) \not= 0$ and thus:

\[(1 - \overline{u}^{-1}\overline{\phi})((1 - \overline{u}^{p^{-n}})^{-1}(\gamma)) = \alpha\]

Finally, setting $\beta = -\overline{u}^{-1}(1 - \overline{u}^{p^{-n}})^{-1}\gamma$ gives us $(\overline{\phi} - \overline{u})(\beta) = \alpha$. As $u \in \mathbb{F}_{p}^{\times}$ then $\beta \in M \subseteq \overline{T}$ and thus we have the surjectivity of $\overline{\phi} - \overline{u}: \overline{T} \to \overline{T}$ in the case where $d = 1$.

\bigskip

\subsection{When  d is greater than 1}

We will now deal with what happens when the dimension of $F$ is any positive integer greater than $1$.

In this case, $u$ is an invertible $d \times d$ matrix over $\mathbb{Z}_{p}$; thus $\overline{u}$ is an invertible $d \times d$ matrix over $\mathbb{F}_{p}$. We shall also write $\phi$ in the map $\phi - u$ as $\Phi$, since it is now a $d \times d$ diagonal matrix with every non-zero entry being the automorphism $\phi$.

If $\overline{u}$ is the identity then the mathematics in chapter $4$ of``Local Fields and their Extensions'' will again showcase why $\overline{\Phi} - \overline{u}$ is surjective \cite{fesenko1}. So, assume $\overline{u} \not= 1$; this means that $\overline{u}$ has eigenvalues not equal to $1$.

By a change of basis we may write $\overline{u}$ as $u' \oplus u''$. Here, $u'$ may be $0$ but is an identity matrix if it has positive dimension. Meanwhile, $u''$ has no eigenvalues equal to $1$. From the above assumptions we know that $u''$ is an invertible matrix of positive dimension $d'' \le d$.

If $A$ is the change of basis matrix that transforms $\overline{u}$ to $u' \oplus u''$ then, as $\overline{u}$ is over $\mathbb{F}_{p}$, $A$ is over $\mathbb{F}_{p}$. $A$ and $A^{-1}$ commute with $\overline{\Phi}$, which is a diagonal matrix with $\overline{\phi}$ at every non-zero value. Therefore, we can use $A$ to get the new map $\overline{\Phi} - (u' \oplus u'')$. Write $\overline{\Phi}$ as $\Phi' \oplus \Phi''$. Here $\Phi'$ has the same dimension as $u'$, and thus may be $0$, and $\Phi''$ has the same dimension as $u''$. This gives us the map $(\Phi' - u') \oplus (\Phi'' - u'')$.

$\Phi' - u'$ has been dealt with previously; so we can focus on $\Phi'' - u''$, which we shall relabel $\Phi - u$. Now, $u$ has no eigenvalues equal to $1$ and thus $u^{p^{n}}$ does not either for all non-negative integers $n$. Therefore, $(1 - u^{p^{n}})$ is an invertible matrix for all non-negative integers $n$. We can now use a nearly identical method as in the case where $d = 1$ to prove the surjectivity of $\Phi - u$.

We are dealing with $(\Phi' - u') \oplus (\Phi'' - u'')$, so we need to change the basis back to what it originally was. Doing so finishes the proof the surjectivity of the map. This gives us that $\overline{\Phi} - \overline{u}$ is surjective for all values of $d$ greater than $1$.

\bigskip

We have now shown that $\overline{\Phi} - \overline{u}$ is surjective for all positive values of $d$ and thus we have proven Lemma $3$.

\bigskip

\section{Finishing Off Theorem 1}

With Lemma $3$ proved we now move onto the square:

\bigskip

\[ \xymatrix{(G^{\mathtt{ab}})\,\widehat{•} \ar[r] \ar[d]^{1 - u_{J}} & (\mathcal{U}^{1}(E)^{d}/I_{G}\mathcal{U}^{1}(E)^{d})\,\widehat{•} \ar[d]^{\Phi - u_{J}} \\
(G^{\mathtt{ab}})\,\widehat{•} \ar[r] & (\mathcal{U}^{1}(E)^{d}/I_{G}\mathcal{U}^{1}(E)^{d})\,\widehat{•}} \]

\bigskip

\noindent which we want to show is commutative. Here the maps in the two rows are derived from the exact sequence created earlier in Lemma $2$. For simplicity we will assume that $d = 1$ as the more general case is much the same.

Now, $E/L$ is unramified; so let $\pi$ be a prime element of $L$, it is also a prime element of $E$. We may use this $\pi$ in the homomorphism between $(G^{\mathtt{ab}})\,\widehat{•} \cong \oplus_{j \in J} (G^{\mathtt{ab}})$ and $(\mathcal{U}^{1}(E)^{d}/I_{G}\mathcal{U}^{1}(E)^{d})\,\widehat{•} \cong \oplus_{j \in J} \mathcal{U}^{1}(E)^{d}/I_{G}\mathcal{U}^{1}(E)^{d}$. This homomorphism sends $\oplus_{j \in J}\alpha_{j} \in (G^{\mathtt{ab}})\,\widehat{•}$ to the element $\oplus_{j \in J} (\pi^{\alpha_{j} - 1} I_{G}\mathcal{U}^{1}(E))$.

Every map of the square acts on a term by term basis in the direct sum so we may pick an arbitrary $j \in J$ and deal only with $\phi_{j}$ and $u_{j}$, which are simplified to $\phi$ and $u$ respectively. We, therefore, want to show that modulo $I_{G}\mathcal{U}^{1}(E)$, we have $(\pi^{\alpha - 1})^{\phi - u} = \pi^{\alpha^{1 - u} - 1}$.

The map from $(G^{\mathtt{ab}})\,\widehat{•}$ to $(\mathcal{U}^{1}(E)^{d}/I_{G}\mathcal{U}^{1}(E)^{d})\,\widehat{•}$ is a homomorphism, so for every $n \in \mathbb{Z}$ we have $(\pi^{\alpha - 1})^{n} = \pi^{\alpha^{n} - 1}$. We also have that $\pi \in L$, which means $\phi(\pi^{\alpha - 1}) = \pi^{\alpha - 1}$ and thus $(\pi^{\alpha - 1})^{\phi - u} = (\pi^{\alpha - 1})^{1 - u}$. We have set $d = 1$ therefore, from our earlier workings, $u$ is an integer such that $1 \le u \le p - 1$. This gives us that $(\pi^{\alpha - 1})^{1 - u} =  \pi^{\alpha^{1 - u} - 1}$, which is precisely what we wanted to show.

To finish off the proof of Theorem $1$ observe the diagram below:

\bigskip

\[ \xymatrix{&& (1) \ar[d] & (1) \ar[d] & \\
&& V(L)/(V(L) \cap I_{G}\mathcal{U}^{1}(E)^{d}) \ar[r]^-{\mathcal{N}}  \ar[d] & V(K) \ar[r] \ar[d] & (1) \\
(1) \ar[r] & (G^{\mathtt{ab}})\,\widehat{•} \ar[r] \ar[d]^{1 - u_{J}} & (\mathcal{U}^{1}(E)^{d}/I_{G}\mathcal{U}^{1}(E)^{d})\,\widehat{•} \ar[d]^{\Phi - u_{J}} \ar[r]^-{\mathcal{N}} & (\mathcal{U}^{1}(T)^{d})\,\widehat{•} \ar[r] \ar[d]^{\Phi - u_{J}} & (1) \\
(1) \ar[r] & (G^{\mathtt{ab}})\,\widehat{•} \ar[r] & (\mathcal{U}^{1}(E)^{d}/I_{G}\mathcal{U}^{1}(E)^{d})\,\widehat{•} \ar[r]^-{\mathcal{N}} \ar[d] & (\mathcal{U}^{1}(T)^{d})\,\widehat{•} \ar[r] \ar[d] & (1) \\
&& (1) & (1) & } \]

\bigskip

Lemmas $2$ and $3$ tell us that the bottom two rows and right two columns are exact. Meanwhile, the work we did earlier this section combined with the fact that $\Phi - u_{J}$ commutes with the elements of $G$, thus with the norm maps labelled $\mathcal{N}$, tells us the diagram is commutative.

We now use the Snake Lemma to get the exact sequence:

\[ \xymatrix{ V(L)/(V(L) \cap I_{G}\mathcal{U}^{1}(E)^{d}) \ar[r]^-{\mathcal{N}} & V(K) \ar[r] & (G^{\mathtt{ab}})\,\widehat{•}\,/((1 - u_{J})(G^{\mathtt{ab}})\,\widehat{•}) \ar[r] & (1) } \]

\noindent This is achieved by taking into account the information shown on the diagram.

This gives the result:

\[ V_{u_{J}}(K)/\mathcal{N}_{L/K}(V_{u_{J}}(L)) \cong ((G^\mathtt{ab})^{d})\,\widehat{•}\,/(I - u_{J})(((G^\mathtt{ab})^{d})\,\widehat{•}\,) \]

\noindent Which finishes to proof of Theorem $1$.

\bigskip

\section{A Local Class Field Theory Aside}

The theorem that we have just proven, and in fact Theorem $1$ in the paper by Lubin and Rosen \cite{Rubin1}, can be looked at from a Local Class Field Theory perspective.

The above theorem sates that there is an isomorphism between $((G^\mathtt{ab})^{d})\,\widehat{•}\,/(I - u_{J})(((G^\mathtt{ab})^{d})\,\widehat{•}\,)$ and $V_{u_{J}}(K)/\mathcal{N}_{L/K}(V_{u_{J}}(L))$. However, it is known that if we let $L/K$ be an abelian finite totally ramified $p$-extension, that Local Class Field Theory, in the case where the residue field of $K$ is only required to be perfect and have characteristic $p > 0$, tells us that there exists a Reciprocity Homomorphism between $(\mathtt{Gal}(L/K)^{\mathtt{ab}})\,\widehat{•}$ and $U_{1, K}/N_{L/K}(U_{1, L})$. This map is in fact an isomorphism \cite{Ivan1}.

The similarity between the two above isomorphisms indicates that the first map can be maybe seen as an analogue of the second. Thus, as the Reciprocity Homomorphism is a fundamental part of Local Class Field Theory, the above Theorem $1$ may be the beginnings of a type of Local Class Field Theory dealing with abelian varieties with good reduction.

We shall call the isomorphism between $((G^\mathtt{ab})^{d})\,\widehat{•}\,/(I - u_{J})(((G^\mathtt{ab})^{d})\,\widehat{•}\,)$ and $V_{u_{J}}(K)/\mathcal{N}_{L/K}(V_{u_{J}}(L))$ the ``Twisted Reciprocity Homomorphism''. This name is taken from section $5$ of chapter $5$ of Professors Fesenko and Vostokov's ``Local Fields and their Extensions'' \cite{fesenko1}.

\bigskip

\section{Getting from Theorem 1 to Theorem 0}

It is quite simple to use Theorem $1$ to construct the exact sequence that we require for Theorem $0$.

Let $A$ be a $d$-dimensional Abelian Variety over $K$ with good ordinary reduction then we know from the subsection on $V(L)$ that $\hat{A}$, where $\hat{A}$ is the formal group of $A$, is toroidal. Taking into account that $A(\mathcal{O}(L)) = A(L)$, here $L/K$ is a totally ramified Galois extension of degree $p^{n}$ for some non-negative integer $n$, we have the exact sequence:

\[ (0) \to \hat{A}(\mathcal{O}(L)) \to A(L) \to A(k) \to (0) \]

\noindent Since $L/K$ is totally ramified, we have $\overline{L} = \overline{K} = k$.

Using the norm map $N_{L/K}$, and noting that the induced map on $A(k)$ is multiplication by $p^{n}$, we get the exact sequence:

\[ \hat{A}(\mathcal{O}(K))/N_{L/K}(\mathcal{O}(L)) \to A(K)/\mathcal{N}_{L/K}(A(L)) \to A(k)/p^{n}A(k) \to (0) \]

From the subsection on $V(L)$, we have a family of twist matrices $u_{J}$ corresponding to $\hat{A}$ such that $\hat{A}(\mathcal{O}(L)) \cong V_{u_{J}}(L) = V(L)$. Using Theorem $1$ and the theory of limits we end up with the exact sequence:

\[ \oplus_{j \in J} \mathbb{Z}_{p}^{d}/(I - u_{J})(\oplus_{j \in J} \mathbb{Z}_{p}^{d}) \to A(K)/\mathcal{N}_{L/K}(A(L)) \to A(k)_{p} \to (0) \]

\noindent as desired.

\bigskip

\section{Quasi-Finite Residue Fields}

There is a type of field, called a quasi-finite field, whose properties, while less restrictive than those of finite fields, do not cover a large amount of perfect fields.

A field $k$ is quasi-finite if it is perfect and if its absolute Galois group is isomorphic to $\widehat{\mathbb{Z}}$. In other words, if every finite extension of $k$ is a Galois cyclic extension and we also have that for every finite cyclic group, $H$, $k$ has a unique extension whose Galois group is isomorphic to $H$ \cite{fesenko1}. If $k$ is a finite field then it is obviously quasi-finite, but $k$ is also quasi-finite if it is a separable, but not necessarily finite, extension of a finite field. Those two cases does not cover all examples of quasi-finite fields.

When working with the mathematics of ordinary Local Class Field Theory, dealing with a Local Field $K$ that has a quasi-finite residue field is not much harder than if we required $\overline{K}$ to be finite. It is, however, much simpler than if $\overline{K}$ is only required to be perfect and have positive characteristic \cite{Ivan1}. This means that it is good to first consider the case where $\overline{K}$ is quasi-finite before dealing with arbitrary perfect residue fields. When dealing with abelian varieties, though, that simplicity in the mathematics fails to emerge.

If we were to restrict $k$ to being quasi-finite the index $J$ would consist of a single element and thus our family of twist matrices, $u_{J}$, and automorphisms, $\phi_{J}$, would become the single twist matrix and automorphism $u$ and $\phi$ respectively. However, we would still be required to use to the mathematics we have here rather than those in Lubin and Rosen's paper. This is because that by assuming $k$ is quasi-finite we are still not requiring $k$ is finite, something that is necessary for the method Lubin and Rosen employ \cite{Rubin1}. This lack of simplicity explains why it is not helpful to use quasi-finite residue fields as a stepping stone before moving onto more general perfect residue fields.

However, as noted, if $k$ is quasi-finite then $J$ consists of a single element which means we may simplify:

\[ \oplus_{j \in J} \mathbb{Z}_{p}^{d}/(I - u_{J})(\oplus_{j \in J} \mathbb{Z}_{p}^{d}) \to A(K)/\mathcal{N}_{L/K}(A(L)) \to A(k)_{p} \to (0) \]

\noindent to get:

\[\mathbb{Z}_{p}^{d}/(I - u)\mathbb{Z}_{p}^{d} \to A(K)/\mathcal{N}_{L/K}(A(L)) \to A(k)_{p} \to (0) \]

\noindent This is the same exact sequence that appears in Lubin and Rosen's paper \cite{Rubin1}. This, since if $k$ is finite then it is quasi-finite, helps show that Theorem $0$ is a generalisation of Mazur's Proposition $4.39$.

\bigskip

\section{Further Directions to Take this Topic}

An astute reader may have noticed that despite me declaring that I was approaching that topic from a Local Class Field Theory perspective I have barely explored the subject in this paper. I have borrowed basic results from ``Local Fields and their Extensions'', but have mainly used the section on Local Class Field Theory for local fields with perfect residue fields as inspiration of how to modify the groups we were looking at. The only Local Class Field Theory result that I have come up with is the ``Twisted Reciprocity Homomorphism'' and even then that was more of an interesting aside that I noted on the way to proving the main result.

This is because, though ``Local Fields and their Extensions'' was useful, Local Class Field Theory is not really necessary to prove the main result of this paper. That is not to say that Local Class Field Theory is superfluous to this topic though, as seen by the twisted reciprocity homomorphism, since we can use it to direct deeper explorations into the topic of abelian varieties over local fields.

The obvious way to start our studies is to look at the twisted reciprocity homomorphism. The regular reciprocity homomorphism has a number of functorial properties that tells us how $U_{1, K}/N_{L/K}(U_{1, L})$ changes when the extension $L/K$ is manipulated, and vice versa \cite{Ivan1}. It would be interesting to see whether these properties hold, or at least something similar, for the twisted reciprocity homomorphism.

For the next line of inquiry we can note that we have, similar to the case where the residue field of $K$ is finite, an Existence Theorem relating abelian Galois extensions of $K$, though only those that are finite totally ramified $p$-extensions, to open subgroups of $U_{1, K}$ \cite{Ivan1}. This relates to what we studied here, where $L/K$ was a totally ramified Galois $p$-extension.

The Existence Theorem states that there are certain open subgroups of $U_{1, K}$ called Normic subgroups. The norm map gives an order reversing bijection between a class of finite abelian totally ramified $p$-extensions of $K$ and Normic subgroups of $U_{1, K}$.

Explicitly, we have that if we fix a prime element $\pi$ of $K$, every Normic subgroup is equal to $N_{L/K}(U_{1, L}) \subseteq U_{1, K}$ for some unique finite abelian totally ramified $p$-extension $L/K$ with $\pi \in N_{L/K}(L^{\times})$. For such an extension the isomorphism $\oplus_{j \in J} \mathtt{Gal}(L/K) \cong U_{1, K}/N_{L/K}(U_{1, L})$, where the index $J$ is the same index that we have seen before, holds. Finally, we have that if $L$ and $L'$ are both finite abelian totally ramified $p$-extensions with $\pi$ contained in both $N_{L/K}(L^{\times})$ and $N_{L'/K}(L'^{\times})$ then $L \subseteq L'$ if and only if $U_{1, L'} \subseteq U_{1, L}$ \cite{Ivan1}.

One might be able to use this as a means to explore $A(K)$, where $A$ is still an abelian variety with good reduction defined over $K$, and to see if there is some class of open subgroups of $A(K)$ such that using the norm map gives us an analogue to the above. We are talking about this being in the sense of there existing a bijection between a definable collection of finite totally ramified $p$-extensions of $K$ and those subgroups. $A(K)/N_{L/K}(A(L))$ being a group in the exact sequence that was the main focus of this paper means exploring the nature of such groups may be something that can expand naturally from the work that we have done here.

\pagebreak

\nocite{*}
\bibliographystyle{plain}
\bibliography{paper_bibliography.bib}

\end{document}